\documentclass[a4paper,UKenglish,cleveref, autoref, thm-restate]{lipics-v2021_hack}

\hideLIPIcs  

\usepackage{amsmath}
\usepackage{cases}
\usepackage{graphicx}
\usepackage{amssymb}
\usepackage{comment}
\usepackage{placeins}
\usepackage{pdfpages}

\usepackage[ruled,vlined,linesnumbered,norelsize]{algorithm2e}
\usepackage{makecell}

\usepackage{comment}
\usepackage{soul}
\usepackage{graphicx} 

\DeclareMathOperator{\ldet}{ldet}

\DeclareMathOperator{\Trace}{Trace}

\makeatletter
\renewcommand*{\top}{%
  {\mathpalette\@transpose{}}%
}
\newcommand*{\@transpose}[2]{%
  \raisebox{\depth}{$\m@th#1\mathsf{T}$}%
}
\makeatother

\title{Branch-and-bound for D-Optimality with fast local search and variable-bound tightening} 

\author{Gabriel Ponte}{Federal University of Rio de Janeiro, Rio de Janeiro, RJ, Brazil}{gabrielponte@poli.ufrj.br}
{https://orcid.org/0000-0002-8878-6647}{supported in part by CNPq GM-GD scholarship 161501/2022-2}

\author{Marcia Fampa}{Federal University of Rio de Janeiro, Rio de Janeiro, RJ, Brazil}{fampa@cos.ufrj.br}
{https://orcid.org/0000-0002-6254-1510}{supported in part by  CNPq grants 305444/2019-0 and 434683/2018-3}
\author{Jon Lee}{University of Michigan, Ann Arbor, MI, USA   }{jonxlee@umich.edu}{https://orcid.org/0000-0002-8190-1091}{supported in part by AFOSR grant FA9550-22-1-0172}

\authorrunning{G. Ponte, M. Fampa, J. Lee} 

\Copyright{Gabriel Ponte, Marcia Fampa and Jon Lee} 

\keywords{D-optimality, local search, branch-and-bound, variable-bound tightening, convex relaxation} 

\category{} 

\relatedversion{} 

\funding{This work is partially based upon work supported by the 
National Science Foundation under Grant No. DMS-1929284 while 
the authors were in residence at the Institute for Computational and Experimental Research in 
Mathematics in Providence, RI, during the Discrete Optimization program.}

\nolinenumbers 

\EventEditors{Loukas Georgiadis }
\EventNoEds{1}
\EventLongTitle{21st Symposium on Experimental Algorithms (SEA 2023)}
\EventShortTitle{SEA 2023}
\EventAcronym{SEA}
\EventYear{2023}
\EventDate{July 24--26, 2023}
\EventLocation{Barcelona, Spain}
\EventLogo{}
\SeriesVolume{XX}
\ArticleNo{Y}

\begin{document}

\maketitle

\begin{abstract}
We apply a branch-and-bound (B\&B) algorithm
to the D-optimality problem based on 
 a convex mixed-integer nonlinear formulation. We discuss possible methodologies to accelerate the convergence of the B\&B algorithm, by combining the use of   different upper  bounds, variable-bound tightening inequalities, and local-search procedures. Different methodologies to compute the determinant of a matrix after a rank-one update are investigated to accelerate the local-searches.  We discuss our findings through  numerical experiments with randomly generated test problem.
\end{abstract}


\section{Introduction} 

We consider the D-Optimality problem formulated as 

\begin{align*}\label{prob}\tag{D-Opt}
\textstyle
&\max \left\{ \ldet\sum_{\ell\in N} x_\ell v_\ell v_\ell^\top \, : \, \mathbf{e}^\top x=s,~  l\leq x \leq u,~ x\in\mathbb{Z}^n
\right\},\\
&\quad=\max \left\{ \ldet \left(
\sum_{\ell\in N} l_\ell v_\ell v_\ell^\top 
+ \sum_{\ell\in N} x_\ell v_\ell v_\ell^\top\right) \, : \, \mathbf{e}^\top x=s-\mathbf{e}^\top l\,,~  0\leq x \leq u-l,~ x\in\mathbb{Z}^n
\right\},
\end{align*}
where $v_\ell \in \mathbb{R}^{m}$, 
for $\ell\in N:=\{1,\ldots,n\}$,
$0\leq l < u\in\mathbb{Z}^n$, with
$\mathbf{e}^\top l \leq s \leq \mathbf{e}^\top u$. 
\ref{prob} is a fundamental problem in statistics, in the area of ``experimental designs'' (see \cite{Puk}, for example).
Defining $A:= (v_1, v_2, \dots, v_n)^\top$, we consider the least-squares regression 
problem $\min_{\theta\in \mathbb{R}^{m}} \|A\theta -y\|_2$, where $y$ is an 
arbitrary response vector. We assume that $A$ has full column rank, and so there is a unique solution to the 
least-squares problem (for each $y$).  But we consider a situation where each $v_\ell$
corresponds to a costly experiment, which could be carried out up to $u_\ell$ times.
Overall, we have a budget to carry out a total of $s(\geq m)$ experiments, and so
we specify the choices by   $x$ (in \ref{prob}). 
For a given feasible solution $\tilde{x}$, we define $A_{\tilde{x}}$
to be a matrix that has $v_\ell^\top$ repeated $\tilde x_\ell$ times, for $\ell\in N$, as its rows.
This leads to the reduced least-squares 
problem $\min_{\theta\in \mathbb{R}^m} \|A_{\tilde x}\theta -y\|_2$.
The generalized variance of the least-squares parameter estimator $\hat \theta$
is inversely proportional to $\det \sum_{\ell\in N}  \tilde{x}_\ell v_\ell v_\ell^\top $
(which is proportional to the volume of a standard ellipsoidal confidence region for $\theta$),
and so \ref{prob} corresponds to picking the set of experiments to 
minimize the generalized variance of the least-squares parameter estimator $\hat \theta$
(see  \cite{Fedorov}, for example). There is a large literature on heuristic 
algorithms for \ref{prob} and its variations.
\cite{Welch} was the first to approach \ref{prob} with an exact branch-and-bound algorithm,
employing a bound based on Hadamard's inequality and another based on continuous relaxation (apparently without using state-of-the art NLP solvers of that time).  
\cite{KoLeeWayne,KoLeeWayne2} proposed a spectral bound
and analytically compared it with the Hadamard bound;
 also see \cite{LeeLind2019}.
\cite{li2022d} applied a local-search procedure and  an exact algorithm to the D-optimal Data Fusion problem, a particular case of the  D-optimality problem where   $\sum_{\ell\in N}l_\ell v_\ell v_\ell^\top $ is positive definite and known as the existing Fisher Information Matrix (FIM). Moreover, the D-optimal Data Fusion problem  consider only the case where the variables are binary, i.e., $l=0$ and $u=\mathbf{e}$.  Although the Data Fusion and the D-optimality problems have similarities, most techniques  used in \cite{li2022d} rely on  the positive definiteness of the existing FIM and cannot be applied to our problem.

Next, we highlight our contributions. We present in this work
\begin{itemize}
\item three local-search heuristics for \ref{prob},
\item five algorithms to construct an initial solution for the local-search procedures,
\item five procedures to compute the determinant of a rank-one update of a given matrix, knowing the determinant of the matrix. These procedures are essential to the successful application of the local-search procedures, 
\item variable-bound tightening (VBT) inequalities, which are constructed  based on a lower bound for \ref{prob} and on the knowledge of a feasible solution for the Lagrangian dual of its continuous relaxation,
\item  a branch-and-bound algorithm based on a convex mixed-integer nonlinear programming formulation  of \ref{prob}. We investigate possible methodologies to accelerate the convergence of the branch-and-bound algorithm, by combining the use of the VBT  inequalities,  local-search procedures, and the use of  the Hadamard and the spectral upper bounds besides  the bound obtained from the continuous relaxation. 
\item numerical experiments with random generated instances where we first compare the use of the different algorithms to compute the determinant of a rank-one update of a matrix inside the local-search procedures. Then, we  compare several versions of the branch-and-bound algorithm  where subsets of the procedures described above are executed.
\end{itemize}

We note that although \cite{Welch} already considered the application of a branch-and-algorithm for D-optimality, the author  did not use variable tightening inequalities based on convex optimization  or investigated 
the linear algebra of doing a fast local search.

A preliminary version of this paper appeared in \cite{PonteFampaLeeSBPO22}.  Here, we suggest two new  algorithms to construct initial solutions to the local-search procedures, we analyse different ways of computing the determinant of a matrix after a rank-one update,  and  we experiment the procedures proposed inside an enhanced branch-and-bound algorithm. 

A similar solution approach has been successfully applied to the related max\-i\-mum-entropy sampling problem (MESP) (see \cite{AFLW_Using,Anstreicher_BQP_entropy,Kurt_linx,FL2022}), where given the covariance matrix $C$ of a  
Gaussian random $n$-vector, one searches for 
a subset of $s$ random variables which maximizes the ``information'' (measured by ``differential entropy'')
(see \cite{SW,CaseltonZidek1984,LeeEnv,FL2022}, for example).

\vspace{0.1in} 
\noindent {\bf Notation.}
 We let $\mathbb{S}^n$  (resp., $\mathbb{S}^n_+$~, $\mathbb{S}^n_{++}$)
 denote the set of symmetric (resp., positive-semidefinite, positive-definite) matrices of order $n$. 
 We let $\mathbf{diag}(x)$ denote the $n\times n$ diagonal matrix with diagonal elements given by the components of $x\in \mathbb{R}^n$.
We denote an all-ones  vector
by $\mathbf{e}$ and an identity matrix by $I$. 
For matrices $A$ and $B$, 
$A\bullet B:=\Trace(A^\top B)$ is the matrix dot-product.
For matrix $A$, we denote row $i$ by $A_{i\cdot}$ and
column $j$ by $A_{\cdot j}$~.

\section{Variable-bound tightening}\label{sec:ineq}

Next, we present a convex continuous relaxation of \ref{prob} and its Lagrangian dual, which will be used for  tightening the bounds on the variables (it may also be used for variable fixing if sufficiently strong), based on general principles of convex MINLP. 
We define  $A:= (v_1, v_2, \dots, v_n)^\top$ and we note that 
$
\textstyle
\sum_{\ell\in N} x_\ell v_\ell v_\ell^\top = A^\top \mathbf{diag}(x) A.
$
Then,   a convex continuous relaxation of \ref{prob} may be formulated as
\begin{equation}\label{cont_rel}
\max \left\{ \ldet \big(A^\top \mathbf{diag}(x) A\big) \, : \, \mathbf{e}^\top x=s, \, l\leq x\leq u, \, x\in \mathbb{R}^n\right\}.
\end{equation}
It is possible to show that the Lagrangian dual  of \eqref{cont_rel} can be formulated as 
\begin{equation}\label{eq:lag_with_theta}
\begin{array}{lll}
&\min &-\ldet \Lambda  + \lambda^\top u  - \theta^\top l + \nu s - {m},\\
&\text{s.t.} 
&\Lambda \bullet v_iv_i^\top - \lambda_i + \theta_i - \nu = 0,\quad i \in N,\\
&&\Lambda \succ 0,\lambda \geq 0, \theta \geq 0.
\end{array}
\end{equation}

In Theorem  \ref{thm:fix_dopt}, we show how to  tighten variables bounds  for \ref{prob} based on  knowledge of a lower bound  and a feasible solution for the dual problem \eqref{eq:lag_with_theta}.

\begin{theorem}\label{thm:fix_dopt}
Let 
\begin{itemize}
    \item LB be the objective-function value of a feasible solution for  \ref{prob};
    \item $(\hat\Lambda,\hat\lambda,\hat\theta,\hat\nu)$ be a feasible solution for  \eqref{eq:lag_with_theta} with objective-function value $\hat\zeta$.
\end{itemize}
Then, for every optimal solution $x^\star$ for \ref{prob}, we have:
\begin{align}
     &x_k^\star \leq l_k + \left\lfloor
  \left(\hat{\zeta}-{LB}\right)/\hat\theta_k
     \right\rfloor ,\quad ~\forall\; k \in N\text{ such that } \hat \theta_k>0,\label{ineq1}\\
     &x_k^\star \geq u_k-  \left\lfloor
    \left(\hat{\zeta}-{LB}\right)/\hat\lambda_k
     \right\rfloor,\quad ~\forall\; k \in N\text{ such that } \hat \lambda_k>0.\label{ineq2}
     \end{align}
\end{theorem}

\section{Local-search heuristics}\label{sec:heur}
We introduce heuristics  to construct a feasible solution to \ref{prob} by applying  a local-search procedure from an initial solution. We propose different ways of constructing the initial solution and  performing the local search. In the next section, we also investigate procedures to update the objective value  inside the local search in order to make it more efficient.  Without loss of generality, we assume that $l=0$ in \ref{prob}. 

\subsection{Initial solutions from the SVD decomposition of $A$}\label{subsec:heur1}

 Next, we  show how we obtain  initial solutions for our local-search  procedures from the real singular-value decomposition (SVD) $A=U\Sigma V^\top$ (see \cite{GVL1996}, for example), where $U\in\mathbb{R}^{n\times n}$, $V\in\mathbb{R}^{m\times m}$ are orthonormal matrices and $\Sigma=\mathrm{diag}(\sigma_1,\sigma_2,\dots,\sigma_m)\in\mathbb{R}^{n\times m}$ ($n\geq m$) with  singular values $\sigma_1\ge\sigma_2\ge\dots\ge\sigma_m\ge0$.
 
  First, to ensure that we start the local-search procedures with a feasible solution for \ref{prob} with finite objective value, we construct a vector $\tilde{x}\in\{0,1\}^n$, such that $\mathbf{e}^\top \tilde{x}=m$ and $A^\top\mbox{diag}(\tilde{x})A\in\mathbb{S}^m_{++}$\,. This is equivalent to choosing $m$ linearly independent rows of $A$, and setting $\tilde{x}$ as the incidence vector for the selected subset of rows. We denote the set of indices of the selected rows by $\tilde{N}$. To select the linearly independent rows,  we use the Matlab function nsub\footnote{\url{www.mathworks.com/matlabcentral/fileexchange/83638-linear-independent-rows-and-columns-generator}} (see \cite{FLPX2021} for details).



 We note that 
 for each  $k \in N$, we have $\sum_{j\in N} U_{jk}^2 = 1$ and $\sum_{j\in N} U_{kj}^2 = 1$. We define
\[
\textstyle
x^0_j := \sum_{k=1}^s U_{jk}^2\,,\quad   j\in N.
\]
We clearly have  $\mathbf{e}^\top x^0 = s$ and  $0 \leq x^0_j \leq 1$, for all $j\in N$. So, $x^0$ is a feasible solution of \eqref{cont_rel}.  

We let $\tau$ be the permutation of the indices in $N$, such that      $x^0_{\tau(1)}\geq x^0_{\tau(2)}\geq\cdots\geq x^0_{\tau(n)}$.
Then, we propose two procedures to construct a feasible solution $\bar{x}$ for \ref{prob}, considering  $\tau$ and $\tilde{x}$.
\begin{itemize}
    \item ``Bin$(x^0)$'':   
    Let $\bar{N}$ be the first $s-m$ indices in $\tau$ (which depends on $x^0$) that are not in $\tilde {N}$. Set  $\bar{x}_{j}:=1$, for $j\in\bar{N}$, and $\bar{x}_{j}:=\tilde{x}_j$\,, for $j\notin \bar{N}$.
     \item ``Int$(x^0)$'': 
     Let $\Delta =u-\tilde{x}$ and $\bar{s}=s-m$. Define, for all $j\in N$, 
     \[
     \textstyle
     \tilde y_{\tau(j)} :=  \min\left\{\Delta_{\tau(j)},\max\left\{0,\bar{s}-\sum_{i=1}^{j-1}\tilde{y}_{\tau(i)}\right\}\right\}.
     \]
     Then, set $\bar{x}:=\tilde{x}+\tilde{y}$.
\end{itemize}
We observe that the objective function of \ref{prob} is given by $\ldet \big(\Sigma^{\top}U^{\top}\mathbf{diag}(x)U\Sigma\big)$, and so  the choice of $x$ is related to the rows of $U\Sigma$. Then, we  also define 
    \[
    \textstyle
    {\hat{x}^0_j := \sum_{i=1}^{m} \big(U_{ji} \Sigma_{ii}\big)^2,\quad   j\in N.}
    \]
    Finally, replacing $x^0$ by $\hat{x}^0$ on the procedures described above, we construct two alternative initial solutions to our local-search procedures. We note that although $\hat{x}^0\geq 0$, it need not be  feasible for \eqref{cont_rel}.

\subsection{Initial solution  from the continuous relaxation}
In Algorithm \ref{alg:lee_cont}, we present how we compute an initial solution to our local-search procedures from a solution to the continuous relaxation \eqref{cont_rel}.
\begin{algorithm}[!ht]
\footnotesize{
\KwIn{a feasible solution $x^{\mathcal{C}}$ to \eqref{cont_rel}}
	\KwOut{a feasible solution $\bar{x}$ to \ref{prob} }
	$\bar{x} := \lfloor x^{\mathcal{C}} \rfloor$\;
	$k := \mathbf{e}^\top \bar{x}$\;
	$x^{f}:= x^{\mathcal{C}} - 
 \bar{x}$\; 
    \While{$k< s$}{
        $\hat\jmath := \mbox{argmax}\{x^f\}$\;
        $\bar x_{\hat \jmath} := \bar x_{\hat \jmath}  + 1$\;
        $x^f_{\hat\jmath}:=0$\;
        $k := k + 1$\;
    }
    \caption{Convert Continuous to Integer}\label{alg:lee_cont}
}
\end{algorithm}
\subsection{Local-search procedures}\label{subsec:heur2}
In Algorithm  \ref{localsearch}, we present the local-search procedures that consider as the criterion for improvement of the given solution, the increase in the  value of the objective function of \ref{prob}. The neighborhood of a given solution $\bar{x}$ is defined by
\[
\mathcal{N}(\bar{x}):= \{y\in\mathbb{Z}^n~:~0\leq y\leq u,~y_i=\bar{x}_i+1,~y_j=\bar{x}_j-1,~ y_k=\bar{x}_k, k\neq i,k\neq j, \forall  i,j\in N\}.
\]
\begin{algorithm}
	\footnotesize{
		\KwIn{ A feasible solution $\bar{x}$ of \ref{prob}}
		\KwOut{ A feasible solution $\bar{x}$ of \ref{prob}, possibly updated }

$x^0:=\bar{x}$\;
$z^0:=\ldet(A^\top \mathbf{diag}(x^0)A)$\;
$flag:=true$\;
\While{$flag$}
{
$flag:=false$\;
\For{ $i=1,\ldots,n$, such that $\bar{x}_i<u_i$}
{
$x_i := \bar{x}_i + 1$\;
\For{ $j=1,\ldots,n$, such that $\bar{x}_j>0$\,, $j\neq i$}
{
$x:=\bar{x}$\;
$x_j := \bar{x}_j - 1$\;
$z:= \ldet(A^\top \mathbf{diag}(x)A) $\label{innerloop}\;
\If {$z>z^0$}
{
$x^0:=x$\;
$z^0:=z$\;
$flag:=true$\;
\If  {``First improvement''}
{
break loops for $i$ and $j$\;
}
}
}
\If  {``First improvement plus'' \& $flag == true$}
{
break loop for $i$\;
}
}
$\bar{x}:=x^0$\;
}
\caption{Local-search procedures \label{localsearch}}
}
\end{algorithm}
We experiment with the three local-search procedures described next.
\begin{itemize}
    \item ``FI'' (Local Search First Improvement): 
    Starting from $\bar{x}$, the procedure visits the solution in $\mathcal{N}(\bar{x})$ with increased objective value with respect to $\bar{x}$, such that $i$ is the least possible index, and $j$ is the least possible index for the given $i$.   
    \item ``FI$^+$'' (Local Search First Improvement Plus): 
       Starting from $\bar{x}$, the procedure visits the solution in $\mathcal{N}(\bar{x})$ with increased objective value with respect to $\bar{x}$, such that $i$ is the least possible index, and $j$ is selected in $N$, as the  index that maximizes the objective value, for the given $i$.   
    \item ``BI'' (Local Search Best Improvement): 
       Starting from $\bar{x}$, the procedure visits the solution in $\mathcal{N}(\bar{x})$ with increased objective value with respect to $\bar{x}$, such that $i$ and $j$ are selected in $N$, as the pair of indices that maximizes the objective value.   
\end{itemize}

\section{Fast local search}\label{sec:fast}

An efficient local search for \ref{prob} is
based on fast computation of $\ldet\left(B+v_iv_i^\top -v_jv_j^\top\right)$,
already knowing  $\ldet B$, where $B:=
\sum_{\ell\in N}  \bar{x}_\ell v_\ell v_\ell^\top$,
for some $\bar x$ that is feasible for \ref{prob} such that 
$\bar{x}+\mathbf{e}_i-\mathbf{e}_j$  is
also feasible for \ref{prob}. 
If $\ldet\left(B+v_iv_i^\top -v_jv_j^\top\right) > \ldet B$,
then $\bar{x}+\mathbf{e}_i-\mathbf{e}_j$ is an improvement on 
$\bar{x}$ in \ref{prob}.


\subsection{Simplest}\label{sec:simplest}
In the simplest algorithm, we form $\hat{B}$ as $B +v_iv_i^\top - v_jv_j^\top$,
and then we calculate the determinant of $\hat B$ in $\mathcal{O}(m^3)$ flops. 


\subsection{Cholesky update}
Let $B=LL^\top$ be the Cholesky factorization of $B$.  The \texttt{lowrankupdate} and \texttt{lowrankdowndate} Julia functions compute the Cholesky factorization of a rank-one update of an $m\times m$ matrix in  $\mathcal{O}(m^2)$ flops. We first apply a \texttt{lowrankupdate}, outside of the 
inner loop of Algorithm \ref{localsearch}, to get a  Cholesky factorization
$\tilde{L}\tilde{L}^\top$ of $B +v_iv_i^\top$ (from the one for $B$).
Then, inside the inner loop, with $i$ fixed
we apply the a \texttt{lowrankdowndate} to get the Cholesky factorization 
$\hat{L}\hat{L}^\top$
of
 $\hat{B}=B +v_iv_i^\top - v_jv_j^\top$ (from the one for
 $B +v_iv_i^\top$). Finally, we have
$\ldet(\hat{B}) = 2 \sum_{\ell=1}^m \log({\hat L}_{\ell\ell})$.
Computing the  Cholesky  factorization of $\hat{B}$ directly
would have instead required 
$\mathcal{O}(m^3)$ flops (see \cite[Sec. 4.2]{GVL1996}).


\subsection{Sherman–Morrison update}
The well-known Sherman-Morrison formula
\[
(M+ab^\top)^{-1} = M^{-1} -\frac{(M^{-1}a)(b^\top M^{-1})}{(1+b^\top M^{-1}a) }
\]
and the well-known  matrix determinant lemma
\[
\det(M+ab^\top) = (1+b^\top M^{-1}a) \det(M)
\]
are useful for rank-one updates of inverses and determinants, respectively,
in $\mathcal{O}(m^2)$ for an order-$m$ matrix.

Outside the inner loop of Algorithm \ref{localsearch},
we can calculate the inverse of $B +v_iv_i^\top$  from the inverse of $B$, using the Sherman-Morrison formula
(setting $M:=B$, $a:=b:=v_i$).
Inside the inner loop (with $i$ fixed), for each $j$ we can calculate $\ldet(\hat B)$ 
from the inverse of   $B +v_iv_i^\top$, using the matrix-determinant lemma
(setting $M:= B +v_iv_i^\top$,~ 
$a:=-v_j$ and $b:=v_j$).

\subsection{SVD Rank-One update}

The next method requires some preprocessing.
For a given feasible solution $\bar{x}$, we define $A_{\bar{x}}\in \mathbb{R}^{s\times m}$
to be a matrix that has $v_\ell^\top$ repeated $\bar x_\ell$ times, for $\ell\in N$, as its rows.
Note that $A_{\bar x}^\top  A_{\bar x} = B$. But rather than working with $B$ directly, we instead work
with $A_{\bar x}$\,. Let $\phi(j)$ be any row index of $A_{\bar x}$
that contains a copy of $v_j^\top$\,. 
Let
\[
    X := A_{\bar x}  + \mathbf{e}_{\phi(j)}\left(v_{i} - v_{j} \right)^\top.
\] 
This rank-1 update of $A_{\bar x}$ is the result of replacing  $v_j^\top$ with 
$v_i^\top$ in row $\phi(j)$ of $A_{\bar x}$\,.
We note that $\ldet(X^\top X) = \ldet (\hat B)$. 


 Let $A_{\bar x} = U\Sigma V^\top$ with $U \in \mathbb{R}^{s \times m}$, $\Sigma \in \mathbb{R}^{m \times m}$ and $V \in \mathbb{R}^{m \times m}$ be the singular value decomposition (SVD) of $A_{\bar x}$. 
 Let $w:=v_i-v_j$\,.
 We are interested in the SVD of
    \begin{equation}\label{eq:upsvdinit}
        A_{\bar x}+\mathbf{e}_{\phi(j)} w^\top = 
    \begin{bmatrix}U&\mathbf{e}_{\phi(j)}\end{bmatrix}
    \begin{bmatrix}
    \Sigma ~&~\mathbf{0}\\
    \mathbf{0}~&~I
    \end{bmatrix}
    \begin{bmatrix}V&w\end{bmatrix}^\top
    \end{equation}
    expressed as modifications to $U,\Sigma,V$.\\
    
From \cite{brand2006fast}, let $p \in \mathbb{R}^{s}$ where $p := (I-UU^\top)\mathbf{e}_{\phi(j)}$, and $K \in \mathbb{R}^{(m + 1) \times m}$, where
     \[
     K:=        \begin{bmatrix}
\Sigma V^\top  + U^\top \mathbf{e}_{\phi(j)} w^\top\\
        \|p\|w^\top
        \end{bmatrix}~.
        \]
   Let
     $K = {\tilde U}{\tilde \Sigma}{\tilde V}^\top$ be the singular value decomposition of $K$.
     Then 
    \begin{equation}\label{eqX}
        X= A_{\bar x} + \mathbf{e}_{\phi(j)} w^\top =  \left(\begin{bmatrix}
        U ~&~ p/\|p\|
        \end{bmatrix}\tilde U\right){\tilde \Sigma} \tilde V^\top~,
    \end{equation}

\noindent and 
\[
\ldet \hat B = \ldet (X^\top X)= 
\ldet\left((A_{\bar x}+\mathbf{e}_{\phi(j)} w^\top)^\top(A_{\bar x}\mathbf{e}_{\phi(j)} w^\top)\right) = \textstyle 2 \sum_{\ell=1}^m \log\left(\tilde \Sigma_{\ell\ell}\right).
\]



Inside the $i,j$  loops in 
Algorithm \ref{localsearch}, working with the
$(m+1)\times m$ matrix $K$,
the $m\times m$ matrix $\Sigma V^\top$ does not change. 
We only need the singular values of $K$ (and not the singular vectors).
In this case, we could employ the 
Golub-Reinsch Algorithm which uses about
 $\frac{8}{3}m^3+4m^2$ flops. 
 
 The direct computation of the singular values of $X$ would, instead, use about $4sm^2-\frac{4}{3}m^3$ if the Golub-Reinsch algorithm was applied (recommended when $s \lessapprox \frac{5}{3}m$), or about $2sm^2+2m^3$ flops if the  R-SVD algorithm was applied (recommended when $s \gtrapprox \frac{5}{3}m$). 

 Outside the $i,j$ loops, we need the complete SVD (singular values and singular vectors) of $X$ to restart the local search.  We can use the SVD of $K$ from the  pair $(i,j)$ that determines the new solution $\bar x$, to make the computation more efficient (see \eqref{eqX}). 
 To compute the complete SVD 
of an $(m+1)\times m$ matrix, we could again employ the 
Golub-Reinsch algorithm which, in this case, uses about $22m^3+14m^2$. 

Computing the complete SVD of an $s\times m$ matrix, would require, instead, about $14sm^2+8m^3$ if the Golub-Reinsch algorithm was applied (recommended when $s \lessapprox \frac{3}{2}m$), or about $6sm^2+20m^3$ flops if the  R-SVD algorithm was applied (recommended when $s \gtrapprox \frac{3}{2}m$) (see \cite[Sec. 5.4.5]{GVL1996}).


Having complexity $\mathcal{O}(m^3)$ per iteration, this algorithm is unlikely to be competitive with the $\mathcal{O}(m^2)$.

\subsection{QR Rank-One update}
Similarly with the preceding section where we updated an SVD factorization, we can compute a QR factorization of $X$ knowing a QR factorization of  $A_{\bar x}$\,. 
Computing the  QR factorization of $X$ directly by the Householder QR algorithm
would have instead required about 
$2m^2(s-m/3)$ flops (see \cite[Sec. 5.2.1]{GVL1996}).
 Let $A_{\bar x} = QR$, where $Q\in\mathbb{R}^{s\times s}$ and $R\in\mathbb{R}^{s\times m}$. The \texttt{qrupdate} Matlab function (which we have implemented in Julia)  computes the QR factorization of a rank-one update of an $s\times m$ matrix  in  $\mathcal{O}(s^2)$ flops, getting $X= \tilde{Q}\tilde{R}$ (see \cite[Sec. 12.5.1]{GVL1996} for the algorithm). Then, we have
 \[
 \ldet\hat B=
 \textstyle 2 \sum_{\ell=1}^m \log\left(\tilde R_{\ell\ell}\right).
 \]
 This algorithm could possibly be competitive with some of the 
  $\mathcal{O}(m^2)$ approaches, when $s$ is not much larger than $m$.

    
    

\section{Hadamard and Spectral Bounds}

Without loss of generality,
we assume that
$l=0$ in \ref{prob}. 
Let $A_u$ denote a $\mathbf{e}^\top u \times m$
matrix obtained by repeating each row $i$ of
$A$ a total of $u_i$ times in $A_u$\,.
In this way, \ref{prob} becomes a 0/1 optimization problem on $A_u$\,, and we can apply some bounds of \cite{KoLeeWayne}. 
The \emph{Hadamard bound}  is defined as
\begin{equation}\label{hada_bound}
\textstyle
\mathcal{H} :=  \sum_{\ell=1}^{s} \log\left(1 +  \phi_{\ell}^2( A_u)\right),
\end{equation}
where $\phi_\ell(A_u)$ denotes  the
denotes the $\ell^{\mbox{th}}$-greatest 2-norm over
the rows of $A_u$\,, and the
\emph{spectral bound} is defined as
\begin{equation}\label{spec_bound}
\textstyle
    \mathcal{S} :=   \sum_{\ell=1}^{s} \log\left(1 +  \sigma_\ell^2( A_u)\right),
\end{equation}
where $\sigma_\ell(A_u)$ denotes the $\ell^{\mbox{th}}$-greatest singular value of  $A_u$\,. \cite{KoLeeWayne}
gives details about how to adapt these bounds,
inside of branch-and-bound
(which becomes complicated when 
the set of rows of $A$ fixed into a solution
do not span $\mathbb{R}^m$). 




\section{Numerical Experiments}\label{sec:num_exp}

For our numerical experiments, we implemented the three local-search procedures described in Subsection \ref{subsec:heur2}, namely, ``FI'', ``FI$^+$'', and ``BI'', and we  initialized each procedure with the four methods proposed in Subsection \ref{subsec:heur1}, namely, ``Bin$(x^0)$'', ``Int$( x^0)$'',``Bin$(\hat{x}^0)$'', ``Int$(\hat{x}^0)$''. All the methods described in Section \ref{sec:fast} to compute the determinant at each iteration of the local-search procedures were implemented and compared to each other. 

We also implemented different versions of a branch-and-bound (B\&B) algorithm to obtain optimal solutions for our test instances, where the bounds are obtained with the convex continuous relaxation \eqref{cont_rel}.   The different versions of the B\&B apply all, some, or none of the  procedures described in the following at each node of the B\&B enumeration tree in an attempt to reduce its size and make the algorithm more efficient. 

\begin{itemize}
    \item VBT: Compute the variable-bound tightening (VBT) inequalities \eqref{ineq1} and \eqref{ineq2} and include them as cuts in the current subproblem, also fixing variables when possible. To compute the inequalities, we use  the values of the optimal dual variables for the continuous relaxation \eqref{cont_rel} of the current subproblem  and the best known lower bound for \ref{prob}.
    \item LSI: Apply the local-search procedures from the solution of the continuous relaxation whenever the solution obtained is integer, in an attempt to increase the lower bound LB on the objective value of \ref{prob}.  
    \item LSC: Apply Algorithm \ref{alg:lee_cont} from the solution of the continuous relaxation whenever the solution is not integer. Then apply the local-search procedures from the integer solution obtained. Unlike we do for LSI, in this case we run the local-search procedures at every node of the B\&B algorithm. 
    \item HS: Compute the Hadamard bound \eqref{hada_bound} and the spectral bound \eqref{spec_bound} and compare them to the upper bound given by the continuous relaxation of the current subproblem, considering the best upper bound when testing if the node can be fathomed. 
\end{itemize}

Algorithm \ref{alg:node_ls}
shows what is executed at each node of the B\&B algorithm when all the enhancement procedures described above are applied.  \begin{algorithm}
	\footnotesize{
\While{$true$}
{
Get $z^\mathcal{H}$ and $z^{\mathcal{S}}$ from the Hadamard and spectral bounds\;
Get ${x},{z}^\mathcal{C},{\lambda},\theta$  from the continuous relaxation \eqref{cont_rel}\;
\If{$x$ is an integer feasible solution}{
$x^0_{LS} := x$\;
break while loop\;
}
\ElseIf {${x}$ is a continuous feasible solution and $\min\{{z}^\mathcal{C},{z}^\mathcal{H},{z}^\mathcal{S}\} > LB$ }{
    Apply VBT\;
    \If{VBT didn't change the bounds of any variable}{
    Get $x_{LS}^0$ from Algorithm \ref{alg:lee_cont} with $x$ as input\;
    break while loop\;
    }
}\Else{
   Node is discarded\;
}
}

Get best $x,z$ from the local-search procedures using $x_{LS}^0$ as input\; 

\If{$z > LB$}{
    $LB := z$\;
}

\caption{Procedure at each node of the enhanced branch-and-bound  \label{alg:node_ls}}
}
\end{algorithm}






The algorithms proposed were coded in Julia v.1.7.1. To solve the convex relaxation \eqref{cont_rel}, we apply Knitro using the Julia package Knitro v0.13.0, and to solve \ref{prob}, we employ  the branch-and-bound algorithm in  Juniper \cite{juniper} (using the \texttt{StrongPseudoCost} branching rule,    and $10^{-5}$ as the tolerance to consider a value as integer). We ran the experiments on
a 16-core machine (running Windows Server 2016 Standard): two Intel Xeon
CPU E5-2667 v4 processors running at 3.20GHz, with 8 cores each, and 128
GB of memory.

To construct our test instances, we used the  Matlab function
\texttt{sprand} to randomly generate 15 $n\times m$ dimensional  matrices $A$ with $m := \lfloor0.25 n\rfloor$ and rank $m$. 
 We generated  three instances for each $n \in \{20,30,50,60,80\}$, we set  $s := 0.5 n$, and used the Matlab function \texttt{randi} to generate a random vector $u$ of each dimension $n$, with integer values  between $1$ and $3$. We set $l=0$ for all instances.    
 
\subsection{Comparing the procedures to update the determinant}

In our first experiments, we verify how the different procedures described in Section \ref{sec:fast}, to update the computation of the determinant at each iteration of the local-search procedures, affect their performance.  We refer the  five  procedures in Section \ref{sec:fast} as ``Simplest'', ``Chol'' (Cholesky), ``SM'' (Sherman-Morrison), ``SVD'', and ``QR''.  

\begin{figure}[!ht]
    \centering
    \includegraphics[scale=0.4]{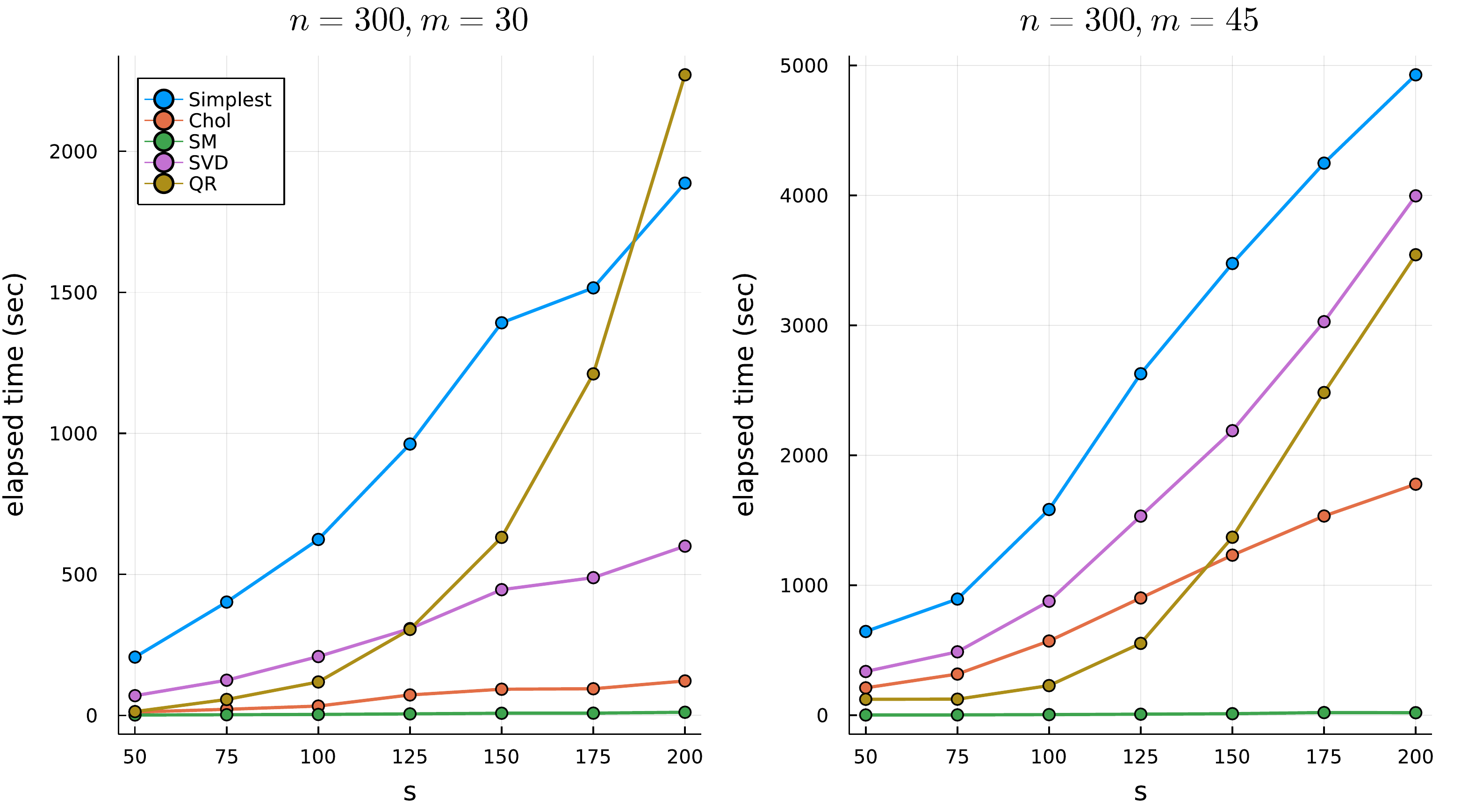}
    \caption{Elapsed times for local-search procedures}
    \label{fig:vary_s_300_30}
\end{figure}

\begin{figure}[!ht]
    \centering
    \includegraphics[scale=0.4]{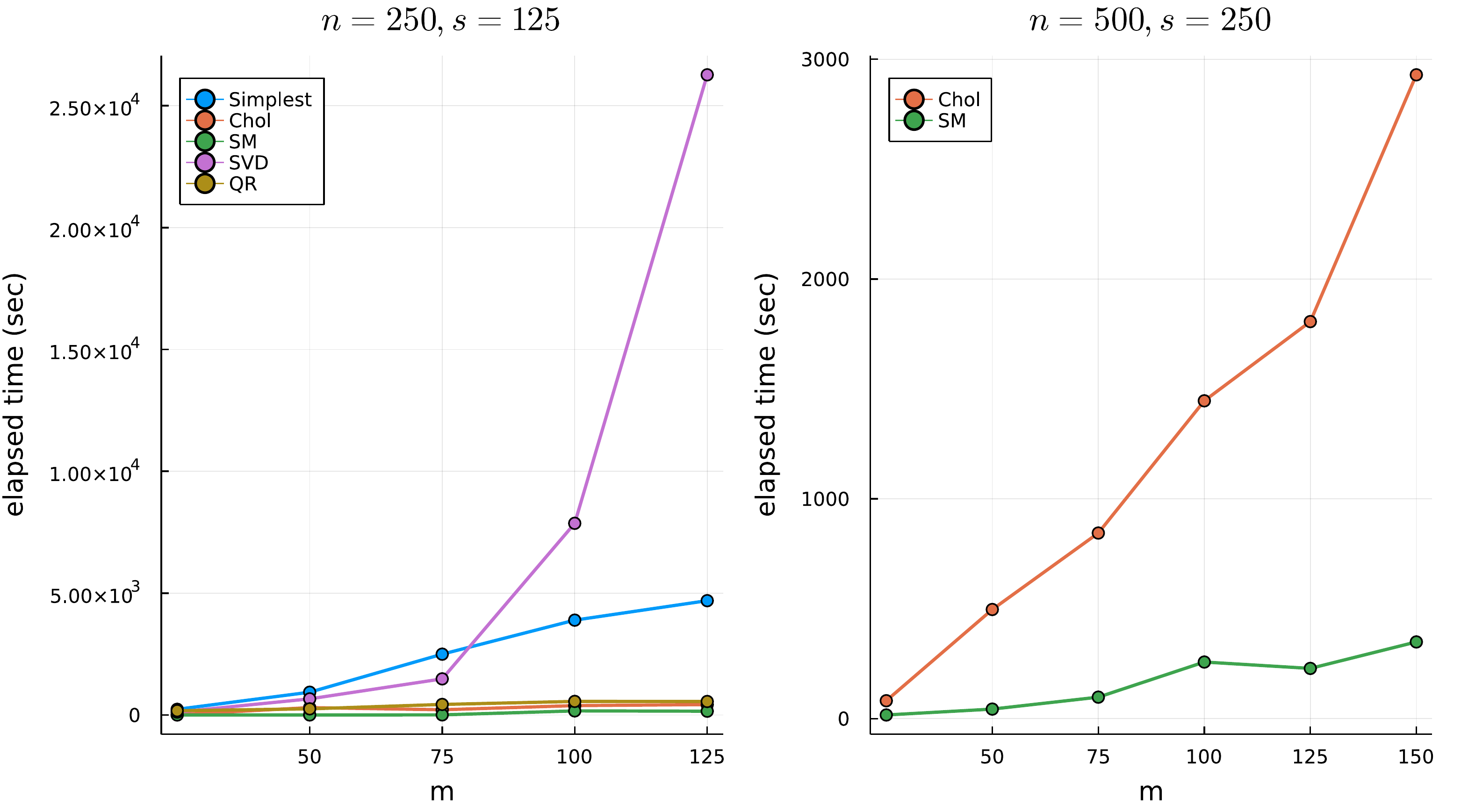}
    \caption{Elapsed times for local-search procedures}
    \label{fig:vary_m_s_125_250}
\end{figure}

\begin{figure}[!ht]
    \centering
    \includegraphics[scale=0.4]{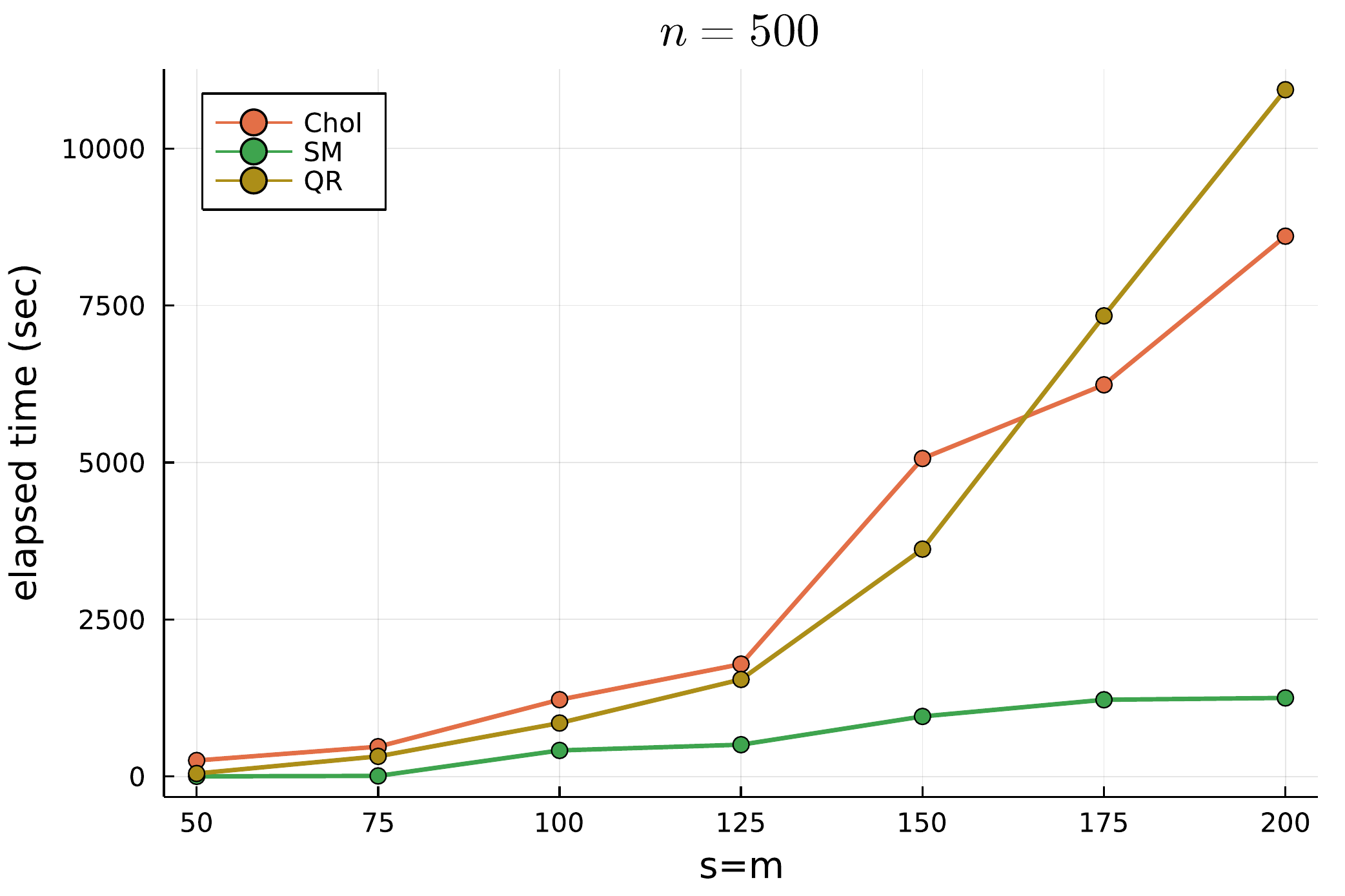}
    \caption{Elapsed times for local-search procedures}
    \label{fig:vary_s_m_n_500_u_1}
\end{figure}

In Figures \ref{fig:vary_s_300_30}, \ref{fig:vary_m_s_125_250} and \ref{fig:vary_s_m_n_500_u_1} we compare the total elapsed times to run the three local-search procedures described in Subsection \ref{subsec:heur2}, starting from the solution ``Bin$(x^0)$'' (see Subsection \ref{subsec:heur1}), using each procedure described in Section \ref{sec:fast}. The times depicted correspond to only one instance generated as described previously, but for these tests we considered other values of $n$, $m$, and $s$, to better observe how the times increase with these parameters. 

From Figure \ref{fig:vary_s_300_30}, we first observe that when $s$ increases in both plots (for $n = 300, m=30$ and $n = 300, m=45$) the times for the QR method (having complexity $\mathcal{O}(s^2)$ per iteration) have a big increase confirming what is expected from theory. We also see that when $m$ increases from $30$ to $45$ all the methods show an increase in time, but for SM we have the smallest times and the smallest increase. It is interesting to note that when $m$ increases from $30$ to $45$, QR becomes more competitive for Chol. In fact, for $m=45$ and $s\leq 125$, QR is faster than Chol, while it is always slower when $m=30$. As QR has complexity $\mathcal{O}(s^2)$ per iteration and Chol has complexity $\mathcal{O}(m^2)$ per iteration,  increasing $m$ is expected to increase more the times for Chol. 
 The results point to SM as the most efficient method in our experiments. 
The second best method can be Chol or QR. As expected, both methods with complexity $\mathcal{O}(m^3)$ per iteration (Simplest and SVD) have bad performance.

The superiority of SM is confirmed in Figure \ref{fig:vary_m_s_125_250}, where we vary $m$. We see that SVD (having complexity $\mathcal{O}(m^3)$ per iteration) becomes very inefficient in comparison to the other methods when $m=100$ or $125$, for 
$n=250, s= 125$. The plot in the right ($n=500, s= 250$) compares the two methods with complexity $\mathcal{O}(m^2)$ per iteration  and we see again a better performance when using SM for this larger instance.

In Figure \ref{fig:vary_s_m_n_500_u_1}, we compare the three methods with complexity $\mathcal{O}(m^2)$ per iteration when $s=m$ (in this case, we set  $u=\mathbf{e}$). In that regime, we confirm that SM is the best option for our use  and that there is no winner between QR and Chol.  

From these experiments, we  see that an efficient implementation of the procedure to update the computation of the determinant can have a significant impact in the efficiency of the local-search procedures.

\subsection{Comparing the versions of branch-and-bound}

In Table \ref{tab:bb}, we analyse the impact of the procedures ``VBT'', ``LSI'', ``LSC'', and ``HS'' on the performance of the branch-and-bound algorithm. On the first row of the table, we identify  the seven  versions of the branch-and-bound algorithm that were executed. On the next four  rows, we identify which procedures are running in the branch-and-bound for each version. Then,  we present the elapsed time to solve the instances  and the  number of nodes on the branch-and-bound enumeration tree  for each version of the algorithm. The parameters $n,m,s$ of each instance are presented in the first column of the table. 

When comparing  versions (1)--(4) to versions (5)--(7), we see that the most successful procedure is LSC. When adding it to the branch-and-bound algorithm both the time and the number of nodes  decrease significantly in general. When applying VBT to the branch-and-bound algorithm already with the LSC procedure included (version (6)), we have another decrease in time and number of nodes for most instances and this is the most successful procedure on our tests. We see that both LCI and HS are not effective in reducing neither the time nor the number of nodes. In fact, we observed that running the local-search procedures only when an integer solution is obtained during the execution of branch-and-bound, rarely improves the current lower bound on the objective value of \ref{prob}. Moreover, the Hadamard and the spectral bounds are rarely stronger than the continuous bound.
\begin{table}[!ht]
\begin{tabular}{c|rrrrrrr}
\hline
BB version&\multicolumn{1}{c|}{(1)}&\multicolumn{1}{c|}{(2)}&\multicolumn{1}{c|}{(3)}&\multicolumn{1}{c|}{(4)}&\multicolumn{1}{c|}{(5)}&\multicolumn{1}{c|}{(6)}&\multicolumn{1}{c}{(7)}\\
\hline
VBT      & \multicolumn{1}{c|}{}        & \multicolumn{1}{c|}{$\checkmark$} & \multicolumn{1}{c|}{$\checkmark$} & \multicolumn{1}{c|}{$\checkmark$} & \multicolumn{1}{c|}{}             & \multicolumn{1}{c|}{$\checkmark$} & \multicolumn{1}{c}{$\checkmark$} \\
LSI      & \multicolumn{1}{c|}{}        & \multicolumn{1}{c|}{}             & \multicolumn{1}{c|}{$\checkmark$} & \multicolumn{1}{c|}{$\checkmark$} & \multicolumn{1}{c|}{}             & \multicolumn{1}{c|}{}             & \multicolumn{1}{c}{}             \\
LSC      & \multicolumn{1}{c|}{}        & \multicolumn{1}{c|}{}             & \multicolumn{1}{c|}{}             & \multicolumn{1}{c|}{}             & \multicolumn{1}{c|}{$\checkmark$} & \multicolumn{1}{c|}{$\checkmark$} & \multicolumn{1}{c}{$\checkmark$} \\
HS       & \multicolumn{1}{c|}{}        & \multicolumn{1}{c|}{}             & \multicolumn{1}{c|}{}             & \multicolumn{1}{c|}{$\checkmark$} & \multicolumn{1}{c|}{}             & \multicolumn{1}{c|}{}             & \multicolumn{1}{c}{$\checkmark$} \\ \hline
$n,m,s$  & \multicolumn{7}{c}{Elapsed time (sec)}                                                                                                                                                                                                              \\ \hline
20,5,10  & \multicolumn{1}{r|}{3.49}    & \multicolumn{1}{r|}{3.50}         & \multicolumn{1}{r|}{3.52}         & \multicolumn{1}{r|}{3.67}         & \multicolumn{1}{r|}{2.68}         & \multicolumn{1}{r|}{2.57}         & 3.18                             \\
20,5,10  & \multicolumn{1}{r|}{0.34}    & \multicolumn{1}{r|}{0.56}         & \multicolumn{1}{r|}{0.45}         & \multicolumn{1}{r|}{0.56}         & \multicolumn{1}{r|}{0.37}         & \multicolumn{1}{r|}{0.40}         & 0.54                             \\
20,5,10  & \multicolumn{1}{r|}{4.73}    & \multicolumn{1}{r|}{5.82}         & \multicolumn{1}{r|}{5.93}         & \multicolumn{1}{r|}{6.16}         & \multicolumn{1}{r|}{5.25}         & \multicolumn{1}{r|}{3.16}         & 3.64                             \\ \hline
30,7,15  & \multicolumn{1}{r|}{161.81}  & \multicolumn{1}{r|}{42.47}        & \multicolumn{1}{r|}{42.74}        & \multicolumn{1}{r|}{45.10}        & \multicolumn{1}{r|}{171.56}       & \multicolumn{1}{r|}{40.72}        & 43.61                            \\
30,7,15  & \multicolumn{1}{r|}{18.90}   & \multicolumn{1}{r|}{11.72}        & \multicolumn{1}{r|}{12.14}        & \multicolumn{1}{r|}{12.59}        & \multicolumn{1}{r|}{22.67}        & \multicolumn{1}{r|}{9.11}         & 11.04                            \\
30,7,15  & \multicolumn{1}{r|}{1.45}    & \multicolumn{1}{r|}{2.47}         & \multicolumn{1}{r|}{2.62}         & \multicolumn{1}{r|}{2.60}         & \multicolumn{1}{r|}{1.73}         & \multicolumn{1}{r|}{1.81}         & 2.66                             \\ \hline
50,12,25 & \multicolumn{1}{r|}{469.31}  & \multicolumn{1}{r|}{353.17}       & \multicolumn{1}{r|}{390.00}       & \multicolumn{1}{r|}{377.00}       & \multicolumn{1}{r|}{467.12}       & \multicolumn{1}{r|}{394.20}       & 389.87                           \\
50,12,25 & \multicolumn{1}{r|}{139.08}  & \multicolumn{1}{r|}{142.24}       & \multicolumn{1}{r|}{153.82}       & \multicolumn{1}{r|}{148.98}       & \multicolumn{1}{r|}{141.37}       & \multicolumn{1}{r|}{145.74}       & 146.67                           \\
50,12,25 & \multicolumn{1}{r|}{34.40}   & \multicolumn{1}{r|}{34.46}        & \multicolumn{1}{r|}{38.22}        & \multicolumn{1}{r|}{36.29}        & \multicolumn{1}{r|}{30.17}        & \multicolumn{1}{r|}{30.57}        & 31.19                            \\ \hline
60,15,30 & \multicolumn{1}{r|}{72.78}   & \multicolumn{1}{r|}{60.10}        & \multicolumn{1}{r|}{64.32}        & \multicolumn{1}{r|}{61.85}        & \multicolumn{1}{r|}{34.91}        & \multicolumn{1}{r|}{51.91}        & 51.66                            \\
60,15,30 & \multicolumn{1}{r|}{27.05}   & \multicolumn{1}{r|}{32.82}        & \multicolumn{1}{r|}{33.79}        & \multicolumn{1}{r|}{33.12}        & \multicolumn{1}{r|}{25.81}        & \multicolumn{1}{r|}{37.64}        & 36.84                            \\
60,15,30 & \multicolumn{1}{r|}{27.69}   & \multicolumn{1}{r|}{37.11}        & \multicolumn{1}{r|}{40.77}        & \multicolumn{1}{r|}{37.23}        & \multicolumn{1}{r|}{35.83}        & \multicolumn{1}{r|}{45.08}        & 45.90                            \\ \hline
80,20,40 & \multicolumn{1}{r|}{421.14}  & \multicolumn{1}{r|}{595.66}       & \multicolumn{1}{r|}{622.47}       & \multicolumn{1}{r|}{591.38}       & \multicolumn{1}{r|}{258.11}       & \multicolumn{1}{r|}{211.35}       & 219.96                           \\
80,20,40 & \multicolumn{1}{r|}{602.43}  & \multicolumn{1}{r|}{1085.66}      & \multicolumn{1}{r|}{1128.50}      & \multicolumn{1}{r|}{1062.40}      & \multicolumn{1}{r|}{320.31}       & \multicolumn{1}{r|}{316.72}       & 318.46                           \\
80,20,40 & \multicolumn{1}{r|}{9579.26} & \multicolumn{1}{r|}{12751.16}     & \multicolumn{1}{r|}{13301.99}     & \multicolumn{1}{r|}{14560.62}     & \multicolumn{1}{r|}{5088.98}      & \multicolumn{1}{r|}{4636.28}      & 4719.10                          \\ \hline
$n,m,s$  & \multicolumn{7}{c}{Number of nodes}                                                                                                                                                                                                                 \\ \hline
20,5,10  & \multicolumn{1}{r|}{191}     & \multicolumn{1}{r|}{125}          & \multicolumn{1}{r|}{125}          & \multicolumn{1}{r|}{125}          & \multicolumn{1}{r|}{121}          & \multicolumn{1}{r|}{103}          & 103                              \\
20,5,10  & \multicolumn{1}{r|}{3}       & \multicolumn{1}{r|}{3}            & \multicolumn{1}{r|}{3}            & \multicolumn{1}{r|}{3}            & \multicolumn{1}{r|}{3}            & \multicolumn{1}{r|}{3}            & 3                                \\
20,5,10  & \multicolumn{1}{r|}{337}     & \multicolumn{1}{r|}{303}          & \multicolumn{1}{r|}{303}          & \multicolumn{1}{r|}{303}          & \multicolumn{1}{r|}{315}          & \multicolumn{1}{r|}{219}          & 219                              \\ \hline
30,7,15  & \multicolumn{1}{r|}{10901}   & \multicolumn{1}{r|}{3113}         & \multicolumn{1}{r|}{3113}         & \multicolumn{1}{r|}{3091}         & \multicolumn{1}{r|}{10901}        & \multicolumn{1}{r|}{3113}         & 3091                             \\
30,7,15  & \multicolumn{1}{r|}{1467}    & \multicolumn{1}{r|}{683}          & \multicolumn{1}{r|}{683}          & \multicolumn{1}{r|}{683}          & \multicolumn{1}{r|}{1467}         & \multicolumn{1}{r|}{683}          & 683                              \\
30,7,15  & \multicolumn{1}{r|}{37}      & \multicolumn{1}{r|}{37}           & \multicolumn{1}{r|}{37}           & \multicolumn{1}{r|}{37}           & \multicolumn{1}{r|}{37}           & \multicolumn{1}{r|}{37}           & 37                               \\ \hline
50,12,25 & \multicolumn{1}{r|}{23245}   & \multicolumn{1}{r|}{17755}        & \multicolumn{1}{r|}{17755}        & \multicolumn{1}{r|}{17755}        & \multicolumn{1}{r|}{23245}        & \multicolumn{1}{r|}{17755}        & 17755                            \\
50,12,25 & \multicolumn{1}{r|}{6027}    & \multicolumn{1}{r|}{5331}         & \multicolumn{1}{r|}{5331}         & \multicolumn{1}{r|}{5331}         & \multicolumn{1}{r|}{5445}         & \multicolumn{1}{r|}{4969}         & 4969                             \\
50,12,25 & \multicolumn{1}{r|}{1513}    & \multicolumn{1}{r|}{1221}         & \multicolumn{1}{r|}{1221}         & \multicolumn{1}{r|}{1221}         & \multicolumn{1}{r|}{1077}         & \multicolumn{1}{r|}{919}          & 919                              \\ \hline
60,15,30 & \multicolumn{1}{r|}{2975}    & \multicolumn{1}{r|}{1557}         & \multicolumn{1}{r|}{1557}         & \multicolumn{1}{r|}{1557}         & \multicolumn{1}{r|}{889}          & \multicolumn{1}{r|}{833}          & 833                              \\
60,15,30 & \multicolumn{1}{r|}{569}     & \multicolumn{1}{r|}{463}          & \multicolumn{1}{r|}{463}          & \multicolumn{1}{r|}{463}          & \multicolumn{1}{r|}{309}          & \multicolumn{1}{r|}{403}          & 403                              \\
60,15,30 & \multicolumn{1}{r|}{531}     & \multicolumn{1}{r|}{591}          & \multicolumn{1}{r|}{591}          & \multicolumn{1}{r|}{591}          & \multicolumn{1}{r|}{531}          & \multicolumn{1}{r|}{591}          & 591                              \\ \hline
80,20,40 & \multicolumn{1}{r|}{6955}    & \multicolumn{1}{r|}{6245}         & \multicolumn{1}{r|}{6245}         & \multicolumn{1}{r|}{6245}         & \multicolumn{1}{r|}{2845}         & \multicolumn{1}{r|}{1897}         & 1897                             \\
80,20,40 & \multicolumn{1}{r|}{9013}    & \multicolumn{1}{r|}{11101}        & \multicolumn{1}{r|}{11101}        & \multicolumn{1}{r|}{11101}        & \multicolumn{1}{r|}{2995}         & \multicolumn{1}{r|}{2857}         & 2857                             \\
80,20,40 & \multicolumn{1}{r|}{174481}  & \multicolumn{1}{r|}{178609}       & \multicolumn{1}{r|}{178609}       & \multicolumn{1}{r|}{178609}       & \multicolumn{1}{r|}{72635}        & \multicolumn{1}{r|}{62905}        & 62905                            \\ \hline
\end{tabular}
\caption{Branch-and-bound (BB) algorithms }\label{tab:bb}
\end{table}

In Table \ref{tab:2}, we present for each instance,  the elapsed time to execute the local-search procedures (``LS''), the  branch-and-bound algorithm, which corresponds to version (6) on Table \ref{tab:bb} (``BB(6)''), and to solve the continuous relaxation \eqref{cont_rel} with Knitro at the root node at the branch-and-bound tree (``${z}^\mathcal{C}$''). We also present the objective value obtained with each procedure. Finally, column ``VBT'' shows the number of VBT inequalities that were effective in tightening the bounds of a variable during the execution of the branch-and-bound algorithm, and column ``FV'' shows the number of variables fixed by the VBT inequalities. The times for the local-search procedures corresponds to the execution of the three procedures described in Section \ref{subsec:heur2}, each starting from each of the four initial solutions presented in Section \ref{subsec:heur1}. The objective value corresponds to the best solution found. We see that the local-search procedures are very fast compared to the branch-and-bound algorithm and obtain solutions of very good quality even for the largest instances. The continuous relaxation is also solved in less than 1 second for all the instances and give tight bounds for the instances tested. Finally, we see that the VBT inequalities are effective and fix a significant number of variables.
\begin{table}[ht]
\small
\begin{tabular}{c|rrr|rrr|rr}
\hline
         & \multicolumn{3}{c|}{Elapsed time (sec)} & \multicolumn{3}{c|}{Objective value} & \multicolumn{1}{c}{VBT}        & \multicolumn{1}{c}{FV}         \\
$n,m,s$  & \multicolumn{1}{c}{LS}         & \multicolumn{1}{c}{BB(6)}           & \multicolumn{1}{c|}{${z}^\mathcal{C}$}      & \multicolumn{1}{c}{LS}         & \multicolumn{1}{c}{BB(6)}        & \multicolumn{1}{c|}{${z}^\mathcal{C}$}      & \multicolumn{2}{c}{BB(6)} \\ \hline
20,5,10  & 0.012 & 2.572    & 0.020 & 5.754  & 5.768  & 5.802  & 71    & 49    \\
20,5,10  & 0.007 & 0.403    & 0.012 & 6.206  & 6.206  & 6.226  & 32    & 32    \\
20,5,10  & 0.007 & 3.155    & 0.016 & 5.645  & 5.696  & 5.735  & 30    & 27    \\ \hline
30,7,15  & 0.041 & 40.715   & 0.027 & 8.484  & 8.484  & 8.549  & 445   & 371   \\
30,7,15  & 0.025 & 9.109    & 0.023 & 8.549  & 8.549  & 8.604  & 192   & 138   \\
30,7,15  & 0.025 & 1.806    & 0.016 & 9.186  & 9.186  & 9.232  & 125   & 89    \\ \hline
50,12,25 & 0.110 & 394.195  & 0.039 & 13.660 & 13.660 & 13.712 & 6310  & 5838  \\
50,12,25 & 0.208 & 145.736  & 0.038 & 13.454 & 13.457 & 13.543 & 2852  & 2493  \\
50,12,25 & 0.110 & 30.568   & 0.042 & 14.100 & 14.104 & 14.180 & 868   & 727   \\ \hline
60,15,30 & 0.341 & 51.914   & 0.059 & 16.744 & 16.750 & 16.779 & 1033  & 834   \\
60,15,30 & 0.424 & 37.639   & 0.070 & 16.960 & 16.975 & 17.017 & 492   & 378   \\
60,15,30 & 0.252 & 45.080   & 0.050 & 17.083 & 17.083 & 17.162 & 591   & 348   \\ \hline
80,20,40 & 0.898 & 211.346  & 0.147 & 21.521 & 21.607 & 21.707 & 3581  & 2696  \\
80,20,40 & 0.881 & 316.715  & 0.147 & 21.411 & 21.553 & 21.670 & 3981  & 3052  \\
80,20,40 & 0.937 & 4636.277 & 0.107 & 21.610 & 21.671 & 21.789 & 55275 & 47219 \\ \hline
\end{tabular}
\caption{Performance of local search, branch-and-bound version (6), and VBT inequalities}\label{tab:2}
\end{table}

\newpage

\section{Conclusion}
Our numerical experiments indicate promising directions to investigate in order to improve the efficiency of the branch-and-bound algorithm to solve  the  D-optimality problem. One possible approach,  for example, is to  introduce the use of bounds from   \cite{li2022d}  when $\sum_{\ell\in \hat{N}} \hat{x}_{\ell}v_{\ell}v_{\ell}^{\top}$ is  positive definite, where $x_\ell$ is fixed at $\hat{x}_{\ell}$ at a given subproblem, for all $\ell\in\hat N \subset N$.  
\bibliography{Dopt_FLP}

\begin{thebibliography}{10}

\bibitem{Anstreicher_BQP_entropy}
Kurt~M. Anstreicher.
\newblock {Maximum-entropy sampling and the Boolean quadric polytope}.
\newblock {\em Journal of Global Optimization}, 72(4):603--618, 2018.

\bibitem{Kurt_linx}
Kurt~M. Anstreicher.
\newblock Efficient solution of maximum-entropy sampling problems.
\newblock {\em Operations Research}, 68(6):1826--1835, 2020.

\bibitem{AFLW_Using}
{Kurt M.} Anstreicher, Marcia Fampa, Jon Lee, and Joy Williams.
\newblock Using continuous nonlinear relaxations to solve constrained
  maximum-entropy sampling problems.
\newblock {\em Mathematical Programming}, 85:221--240, 1999.

\bibitem{brand2006fast}
Matthew Brand.
\newblock Fast low-rank modifications of the thin singular value decomposition.
\newblock {\em Linear Algebra and its Applications}, 415(1):20--30, 2006.

\bibitem{CaseltonZidek1984}
William~F. Caselton and James~V. Zidek.
\newblock Optimal monitoring network design.
\newblock {\em Statistics and Probability Letters}, 2:223--227, 1984.

\bibitem{FL2022}
Marcia Fampa and Jon Lee.
\newblock {\em Maximum-Entropy Sampling: Algorithms and Application}.
\newblock Springer, 2022.

\bibitem{FLPX2021}
Marcia Fampa, Jon Lee, Gabriel Ponte, and Luze Xu.
\newblock Experimental analysis of local searches for sparse reflexive
  generalized inverses.
\newblock {\em Journal of Global Optimization}, 81:1057--1093, 2021.

\bibitem{Fedorov}
Valerii~V. Fedorov.
\newblock {\em Theory of optimal experiments}.
\newblock Academic Press, New York-London, 1972.
\newblock Translated from the Russian and edited by W. J. Studden and E. M.
  Klimko.

\bibitem{GVL1996}
Gene~H. Golub and Charles~F. Van~Loan.
\newblock {\em Matrix Computations (3rd Ed.)}.
\newblock Johns Hopkins University Press, Baltimore, MD, USA, 1996.

\bibitem{KoLeeWayne2}
Chun-Wa Ko, Jon Lee, and Kevin Wayne.
\newblock A spectral bound for {D}-optimality, 1994.
\newblock Unpublished.

\bibitem{KoLeeWayne}
Chun-Wa Ko, Jon Lee, and Kevin Wayne.
\newblock Comparison of spectral and {H}adamard bounds for {D}-optimality.
\newblock In {\em M{ODA} 5}, Contrib. Statist., pages 21--29. Physica,
  Heidelberg, 1998.

\bibitem{juniper}
Ole Kröger, Carleton Coffrin, Hassan Hijazi, and Harsha Nagarajan.
\newblock Juniper: An open-source nonlinear branch-and-bound solver in julia.
\newblock In {\em Integration of Constraint Programming, Artificial
  Intelligence, and Operations Research}, pages 377--386. Springer
  International Publishing, 2018.

\bibitem{LeeEnv}
Jon Lee.
\newblock Maximum entropy sampling.
\newblock In A.H. El-Shaarawi and W.W. Piegorsch, editors, {\em Encyclopedia of
  Environmetrics, 2nd ed.}, pages 1570--1574. Wiley, Boston, 2012.

\bibitem{LeeLind2019}
Jon Lee and Joy Lind.
\newblock Generalized maximum-entropy sampling.
\newblock {\em INFOR: Information Systems and Operational Research},
  58(2):168--181, 2020.

\bibitem{li2022d}
Yongchun Li, Marcia Fampa, Jon Lee, Feng Qiu, Weijun Xie, and Rui Yao.
\newblock D-optimal data fusion: Exact and approximation algorithms, 2022.
\newblock Preprint arXiv:2208.03589.

\bibitem{PonteFampaLeeSBPO22}
Gabriel Ponte, Marcia Fampa, and Jon Lee.
\newblock Exact and heuristic solution approaches for the {D}-optimality
  problem.
\newblock In {\em Proceedings of the LIV Brazilian Symposium on Operations
  Research}, volume~54. SOBRAPO, Rio de Janeiro, RJ, Brazil, Nov 2022.

\bibitem{Puk}
Friedrich Pukelsheim.
\newblock {\em Optimal Design of Experiments}, volume~50 of {\em Classics in
  Applied Mathematics}.
\newblock Society for Industrial and Applied Mathematics (SIAM), Philadelphia,
  PA, 2006.
\newblock Reprint of the 1993 original.

\bibitem{SW}
Michael~C. Shewry and Henry~P. Wynn.
\newblock Maximum entropy sampling.
\newblock {\em Journal of Applied Statistics}, 46:165--170, 1987.

\bibitem{Welch}
William~J. Welch.
\newblock Branch-and-bound search for experimental designs based on
  {D}-optimality and other criteria.
\newblock {\em Technometrics}, 24(1):41--48, 1982.

\end{thebibliography}

\end{document}